\documentclass[titlepage,twoside,12pt]{article}
\usepackage{amssymb}
\usepackage{amsfonts}
\begin{document}

{\Large \bf Nel's category theory based \\ \\ differential and integral Calculus, \\ \\
\hspace*{4cm} or \\ \\
Did Newton know category theory ?} \\ \\

Elem\'{e}r E Rosinger \\
Department of Mathematics \\
and Applied Mathematics \\
University of Pretoria \\
Pretoria \\
0002 South Africa \\
eerosinger@hormail.com \\ \\

{\bf Abstract} \\

In a series of publications in the early 1990s, L D Nel set up a study of non-normable
topological vector spaces based on methods in category theory. One of the important results
showed that the classical operations of derivative and integral in Calculus can in fact be
obtained by a rather simple construction in categories. Here we present this result in a
concise form. It is important to note that the respective differentiation does {\it not} lead
to any so called generalized derivatives, for instance, in the sense of distributions,
hyperfunctions, etc., but it simply corresponds to the classical one in Calculus. \\
Based on that categorial construction, Nel set up an infinite dimensional calculus which can
be applied to functions defined on non-convex domains with empty interior, a situation of
great importance in the solution of partial differential equations. \\ \\

{\large \bf 1. The Setup} \\

The presentation follows mainly Nel [1], where further details as well as the proofs can be
found. For convenience, a few basic concepts needed from Category Theory are recalled in short
in the Appendix. \\

Let $I \subset \mathbb{R}$ be any compact interval and let $E$ be a real Banach space. We
consider the Banach space structure on the set ${\cal C} ( I, E )$ of all continuous functions
$f : I \longrightarrow E$, induced by the ${\cal L}^\infty$ norm $||~||_\infty$. \\

Our first aim is to define a {\it derivative} for functions $f$ in a suitable subset of
${\cal C} ( I, E )$, and do so by means of a simple construction in categories. Among such
functions $f$ for which a derivative can be defined by categorial means are those called {\it
paths} in $E$, and they are characterized as follows. \\

Let us consider in ${\cal C} ( I, E )$ the subset \\

(1.1) ~~~~ $ Path~{\cal C} ( I, E ) $ \\

of {\it paths} in $E$, given by all $f \in {\cal C} ( I, E )$, for which there exists $h_f \in
{\cal C} ( I \times I, E )$, such that \\

(1.2) ~~~~ $ f(y) - f(x) ~=~ (y-x)h_f(x,y),~~~ x,~ y \in I $ \\

We give several useful properties of such paths. For that purpose, similar with ${\cal C}
( I, E )$, we first define a Banach space structure on the set ${\cal C} ( I \times I, E )$ of
all continuous functions $h : I \times I \longrightarrow E$. \\

Let us now consider in ${\cal C} ( I \times I, E )$ the following {\it closed} subspace \\

(1.3) ~~~~ $ ad~{\cal C} ( I \times I, E ) $ \\

given by all $h \in {\cal C} ( I \times I, E )$, such that \\

(1.4) ~~~~ $ (y-x)h(x,y) + (z-y)h(y,z) + (x-z)h(z,x) ~=~ 0,~~~ x,~ y,~ z \in I $ \\

It is easy to see that, with the notation in (1.2), we have \\

(1.5) ~~~~ $ \begin{array}{l}
                \forall~~~ f \in path~{\cal C} ( I, E ) ~: \\ \\
                ~~~~~~~ *)~~~ h_f \in ad~{\cal C} ( I \times I, E ) \\ \\
                ~~~~ **)~~~ h_f ~~\mbox{is {\it unique}}
             \end{array} $ \\

Furthermore, if $f \in path~{\cal C} ( I, E )$, then the classical derivative of Calculus,
namely, $f^\prime ( x )$ exists for every $x \in I$, and \\

(1.6) ~~~~ $ \begin{array}{l}
                    f^\prime ( x ) ~=~ h_f ( x, x ) ~=~ \\
                 ~~~~=~ \lim_{~y \to x,~z \to x,~ y \neq z}~ ( f ( z ) - f ( y ) ) /
                                    ( z - y ),~~~ x \in I
              \end{array} $ \\

This property, however, will {\it not} be used in the categorial definition of the derivative
given in the sequel, and instead, it will follow from that definition. \\ \\

{\bf Remark} \\

1) Regardless of (1.6), in the {\it purely algebraic} condition (1.4) which defines
the elements in $ad~{\cal C} ( I \times I, E )$, there is {\it no} involvement of any kind of
operation of limit, derivation, or even merely of a division with a quantity which may
eventually vanish. And as we shall see in the sequel, that purely algebraic condition (1.4) is
{\it sufficient} in order to define and effectively compute both the derivative and the
integral by using a simple {\it diagonal} construction in a suitable category. \\

2) However, the idea in the rather unusual definition (1.3), (1.4) is, in view of (1.2), (1.6),
quite straightforward. Namely, for every function $f : I \longrightarrow E$, we can define the
function \\

(1.7) ~~~~ $ h : I \times I \longrightarrow E $ \\

by the expression \\

(1.8) ~~~~ $ h ( x, y ) ~=~ \begin{array}{|r} ~~
            ( f ( y ) - f ( x ) ) / ( y - x )~~~ \mbox{if}~~ x,~ y \in I,~~ x \neq y \\ \\
            \mbox{an arbitrary element}~ a_x \in E~~~ \mbox{if}~~ x = y \in I
                               \end{array} $ \\

and then, clearly, $h$ satisfies (1.2) and (1.4). \\

The issue, however, is whether $h$ is continuous on $I \times I$, more precisely, whether in
(1.8) one can choose $a_x \in E$, with $x \in I$, so that $h$ becomes continuous. And
obviously, this can be done, if and only if $h \in ad~{\cal C} ( I \times I, E )$, in which
case $f$ is a path and $h = h_f$, thus (1.6) holds.

\hfill $\Box$ \\

Let us note two further useful properties of the functions in \\ $ad~{\cal C}
( I \times I, E )$. \\

First, it is easy to see that every $h \in ad~{\cal C} ( I \times I, E )$ is {\it
symmetric}. \\

Second, given two compact intervals $I=[a,b],~ J=[b,c] \subset \mathbb{R}$ together with
functions $h \in ad~{\cal C} ( I \times I, E ),~ k \in ad~{\cal C} ( J \times J, E )$. If
$h(b,b)=k(b,b)$, then there exists a {\it unique} function $l \in ad~{\cal C} ( ( I \cup J )
\times ( I \cup J ), E )$, such that $l_{I \times I} = h$ and $l_{J \times J} = k$. \\ \\

{\large \bf 2. A Few Related Categories} \\

Let us denote by ${\cal B}an$ the category of Banach spaces and their continuous linear
mappings. \\

For every compact interval $I \subset \mathbb{R}$, we obtain two {\it covariant functors} \\

(2.1) ~~~~ $ \begin{array}{l}
                   {\cal C} ( I, - ) : {\cal B}an \longrightarrow {\cal B}an \\ \\
                   {\cal C} ( I \times I, - ) : {\cal B}an \longrightarrow {\cal B}an
              \end{array} $ \\

acting as follows. Given two Banach spaces $E$ and $F$ and a continuous linear mapping $u :
E \longrightarrow F$, then \\

(2.2) ~~~~ $ \begin{array}{l}
                   {\cal C} ( I, u ) : {\cal C} ( I, E ) \ni f \longmapsto u \circ f \in
                                                    {\cal C} ( I, F ) \\ \\
                   {\cal C} ( I \times I, u ) : {\cal C} ( I \times I, E ) \ni h \longmapsto
                                        u \circ h \in {\cal C} ( I \times I, F )
              \end{array} $ \\

Consequently, in view of (1.1), we can define the {\it covariant functor} \\

(2.3) ~~~~ $  ad~{\cal C} ( I \times I, - ) : {\cal B}an \longrightarrow {\cal B}an $ \\

as the restriction of the functor ${\cal C} ( I \times I, - ) : {\cal B}an \longrightarrow
{\cal B}an$, namely, for $u : E \longrightarrow F$ as above, we have \\

(2.4) ~~~~ $ ad~{\cal C} ( I \times I, u ) : ad~{\cal C} ( I \times I, E ) \ni h \longmapsto
                                             u \circ h \in ad~{\cal C} ( I \times I, F ) $ \\

Finally, let us define the {\it diagonal evaluation} mapping \\

(2.5) ~~~~ $ ed_{~I,~E} : ad~{\cal C} ( I \times I, E ) \longrightarrow {\cal C} ( I, E ) $ \\

by \\

(2.6) ~~~~ $ ( ed_{~I,~E}~ ( h ) ) ( x ) ~=~ h ( x, x ) $ \\

for $h \in ad~{\cal C} ( I \times I, E )$ and $x \in I$. \\

The basic result is given in \\

{\bf Theorem on Diagonal Evaluation (L D Nel)} \\

For every compact interval $I \subset \mathbb{R}$, the mapping $ed_{I,-}$ is a {\it natural
isomorphism} and {\it isometry} between the functors $ad~{\cal C} ( I \times I, - )$ and
${\cal C} ( I, - )$. \\ \\

{\large \bf 3. The Derivative} \\

Given a compact interval $I \subset \mathbb{R}$ and a Banach space $E$, the {\it derivative}
$D$ is defined as the mapping \\

(3.1) ~~~~ $ D : Path~{\cal C} ( I, E ) \longrightarrow {\cal C} ( I, E ) $ \\

where for $f \in Path~{\cal C} ( I, E )$ we have, see (1.5) \\

(3.2) ~~~~ $ D f ~=~ ed_{~I,~E}~ ( h_f ) $ \\

Here the {\it important} point to note is that, as mentioned at pct. 1 in the Remark above,
this definition of derivative does {\it not} refer in any way to limits, derivation, or even
merely to division with a quantity which may eventually vanish. Yet, as seen next, this
concept of derivative {\it recovers} the usual derivative in Calculus for functions $f : I
\longrightarrow E$, see Nel [1] for further details. \\

Indeed, in view of (2.6) and (1.6), we have \\

(3.3) ~~~~ $ ( D f ) ( x ) ~=~ h_f ( x, x ) ~=~ f^\prime ( x ),~~~ x \in I $ \\

Therefore it follows that for every path $f \in Path~{\cal C} ( I, E )$ and $x \in I$, we
have \\

(3.4) ~~~~ $ ( D f ) ( x ) ~=~ \lim_{~y \to x,~z \to x,~ y \neq z}~ ( f ( z ) - f ( y ) ) /
                                    ( z - y ) $ \\

which is the classical definition of derivative in Calculus. \\

Here we give a simple example which shows that the derivative in (3.1), (3.2) does {\it not}
lead to any concept of derivation is a generalized sense, like for instance, distributional.
Instead, it corresponds to the classical derivative in Calculus. \\

Let $I = [ - 1, 1]$, $E = \mathbb{R}$, and let us take the function $f \in {\cal C} ( I, E )$
defined by $f = x_{+}$, that is, $f(x)=0$, for $x \in [-1,0]$, while $f(x)=x$, for $x \in
[0,1]$. \\

Then $f \in {\cal C} ( I, E )$, however $f$ is {\it not} a path, therefore, the derivation in
(3.1) - (3.3) {\it cannot} be applied to it. \\

Indeed, assume that $f$ is a path. Then (1.4) gives $h_f \in ad~{\cal C} ( I \times I, E )$
such that $f(y)-f(x)=(y-x)h_f(x,y)$, for $x,~y \in I$. Consequently (3.4), (3.5) imply that

$$ \lim_{~x \to 0,~x \neq 0}~ f(x) / x ~=~ h_f(0,0) $$

which is absurd, since the limit in the left hand term does not exist. \\

It follows that the that categorial definition of derivative in (3.1), (3.2) rather corresponds
to the usual derivative of ${\cal C}^1$-smooth functions. \\

Let us further illustrate the above categorial concept of derivative in (3.1), (3.2) in the
case of some simple functions. \\

A function $f \in {\cal C} ( I, E )$ is called {\it affine}, if and only if, for suitable
$c, d \in E$, we have \\

(3.5) ~~~ $ f ( x ) = c + x d,~~~ x \in I $ \\

We denote by $Aff~( I, E )$ the subset of all affine functions $f \in {\cal C} ( I, E )$. \\

A function $f \in {\cal C} ( I, E )$ is called {\it polygonal}, if and only if it is
piece-wise affine on $I$. We denote by $Poly~( I, E )$ the subset of all polygonal functions
$f \in {\cal C} ( I, E )$. \\

It is well known, Dieudonne, that the set of polygonal functions $Poly~( I, E )$ is {\it
dense} in ${\cal C}( I, E )$. \\

Now for every affine function $f \in Aff~( I, E )$ let us define the corresponding {\it
averaging} function \\

(3.6) ~~~ $ av_f \in {\cal C} ( I \times I, E ) $ \\

by \\

(3.7) ~~~ $ av_f ( x, y ) ~=~ c + ( ( x + y ) / 2 ) d,~~~ x, y \in I $ \\

Then obviously \\

(3.8) ~~~ $ av_f \in ad~{\cal C} ( I \times I, E ),~~~ ed_{~I,~E}~ ( av_f ) ~=~ f $ \\

It follows therefore, Nel [1, p.53], that one can in a piece-wise manner extend the averaging
function in (3.6), (3.7) to all polygonal functions $f \in Poly~( I, E )$, and obtain
relations similar to those in (3.8). \\

Moreover, since the polygonal functions are dense in ${\cal C} ( I, E )$, it follows that one
can further extend the averaging function in (3.6), (3.7) to the whole of ${\cal C} ( I, E )$,
with the preservation of the relations in (3.8). \\ \\

{\large \bf 4. The Integral} \\

Based on the categorial concept of derivative in (3.1), (3.2), we can define as well a notion
of {\it integral}. Indeed, such an integral is given by the mapping \\

(4.1) ~~~ $ \int : I \times I \times {\cal C} ( I, E ) \ni ( a, b, f ) ~\longmapsto~
                          \int_a^b f ~\stackrel{def}=~ ( b - a )\, av_f ( a, b ) \in E $ \\

Now, the following two properties result : \\

For every path $f \in Path~( I, E )$ and every $a, b \in I$, we have the classical
Newton-Leibniz formula of Calculus \\

(4.2) ~~~ $ \int^b_a f^\prime ~=~ f ( b ) - f ( a ) $ \\

where the derivative $f^\prime$ is given by the categorial concept in (3.1), (3.2). \\

Conversely, for every given $f \in {\cal C} ( I, E )$ and $a \in I$, let us define $F \in
{\cal C} ( I, E )$ by \\

(4.3) ~~~ $ F ( x ) ~=~ \int^x_a f,~~~ x \in I $ \\

Then $F$ is a path in ${\cal C} ( I, E )$, that is, $F \in Path~( I, E )$, and with the
categorial concept of derivative in (3.1), (3.2), we have \\

(4.4) ~~~ $ F^\prime ~=~ f $ \\ \\

{\large \bf Appendix} \\

One way to define a {\it category} ${\cal A}$ is as being a directed graph which has the
following properties : \\

C1)~ The vertices are called {\it objects}, while the directed edges are called {\it morphisms},
and given two objects $A$ and $B$, the totality of morphisms $A \stackrel{f}{\longrightarrow}
B$ is a set denoted by $Hom_{\cal A}(A,B)$. \\

C2)~ For every object $A$ there exists a morphism $A \stackrel{id_A}{\longrightarrow} A$ with
property C3.2) below. \\

C3)~ For every two morphisms $A \stackrel{f}{\longrightarrow} B \stackrel{g}{\longrightarrow}
C$, there exists a well defined morphism $A \stackrel{g \circ f }{\longrightarrow} C$, such
that : \\

\hspace*{0.5cm} C3.1)~ the operation $\circ$ is associative \\

\hspace*{0.5cm} C3.2)~ for every morphism $A \stackrel{f}{\longrightarrow} B$, we have $f \circ
id_A = f = id_B \circ f$ \\ \\

A {\it functor} $F$ acts between two categories ${\cal A}$ and ${\cal X}$, namely \\

$~~~~~~~~~~~~~~~~~~~~~~~~ {\cal A} \stackrel{F}{\longrightarrow} {\cal X} $ \\

as follows :

\bigskip
\begin{math}
\setlength{\unitlength}{0.2cm}
\thicklines
\begin{picture}(90,32)

\put(5,29){$A$}
\put(25,26){$f$}
\put(8,29.5){\vector(1,0){36.5}}
\put(46,29){$B$}
\put(5.5,27){\vector(0,-1){20}}
\put(3,16){$F$}
\put(46.5,27){\vector(0,-1){20}}
\put(48,16){$F$}
\put(2,4){$X = F ( A )$}
\put(13,4.5){\vector(1,0){29}}
\put(26.5,16){$F$}
\put(23,7){$\xi = F ( f )$}
\put(44,4){$Y = F ( B )$}
\put(25.5,24){\vector(0,-1){14.5}}

\put(60,29){${\cal A}$}
\put(60.5,27){\vector(0,-1){20}}
\put(60,4){${\cal X}$}
\put(62,16){$F$}

\end{picture}
\end{math}

That is, $F$ takes objects $A$ in category ${\cal A}$ into objects $X = F ( A )$ in category
${\cal X}$. \\
Further, $F$ takes morphisms $A \stackrel{f}{\longrightarrow} B$ in category ${\cal A}$ into
morphisms $X = F ( A ) \stackrel{\xi = F ( f )}{\longrightarrow} Y = F ( B )$ in category
${\cal X}$. \\

The above actions of the functor $F$ have the properties : \\

For every object $A$ in category ${\cal A}$, we have in category ${\cal X}$ the relation \\

F1)~ $ F ( id_A ) ~=~ id_{F ( A )} $ \\

Also, for every two morphisms $A \stackrel{f}{\longrightarrow} B \stackrel{g}
{\longrightarrow} C$ in category ${\cal A}$, we have in category ${\cal X}$ the relation \\

F2a)~ $ F ( g \circ f ) ~=~ F ( g ) \circ F ( f ) $ \\

in which case the functor $F$ is called {\it covariant}. \\

If on the other hand, we have in category ${\cal X}$ the relation \\

F2b)~ $ F ( g \circ f ) ~=~ F ( f ) \circ F ( g ) $ \\

then the functor $F$ is called {\it contravariant}. \\

A {\it natural transformation} $~\nu~$ acts between two functors $F$ and $G$ which, on their
turn, act between the same two categories ${\cal A}$ and ${\cal X}$, namely

\bigskip
\begin{math}
\setlength{\unitlength}{0.2cm}
\thicklines
\begin{picture}(90,22)

\put(10,20){${\cal A}$}
\put(10.5,18){\vector(0,-1){14}}
\put(12,10.5){$F$}
\put(10,1){${\cal X}$}
\put(15,11.5){\vector(1,0){20}}
\put(24,13){$\nu$}
\put(37,10.5){$G$}
\put(40,20){${\cal A}$}
\put(40.5,18){\vector(0,-1){14}}
\put(40,1){${\cal X}$}

\end{picture}
\end{math} \\

And the action of $\nu$ is as follows : \\

N1)~ To every object $A$ in category ${\cal A}$ corresponds in category ${\cal X}$ a
morphism \\

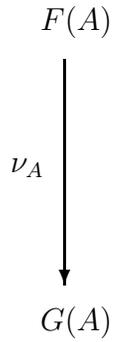
\begin{math}
\setlength{\unitlength}{0.2cm}
\thicklines
\begin{picture}(90,21)

\put(8,20){$F(A)$}
\put(9.6,18){\vector(0,-1){15}}
\put(6,10.5){$\nu_A$}
\put(8,0){$G(A)$}

\end{picture}
\end{math} \\

N2)~ Every morphisms $A \stackrel{f}{\longrightarrow} B$ in category ${\cal A}$ generates
the {\it commutative diagram} in category ${\cal X}$, given by \\

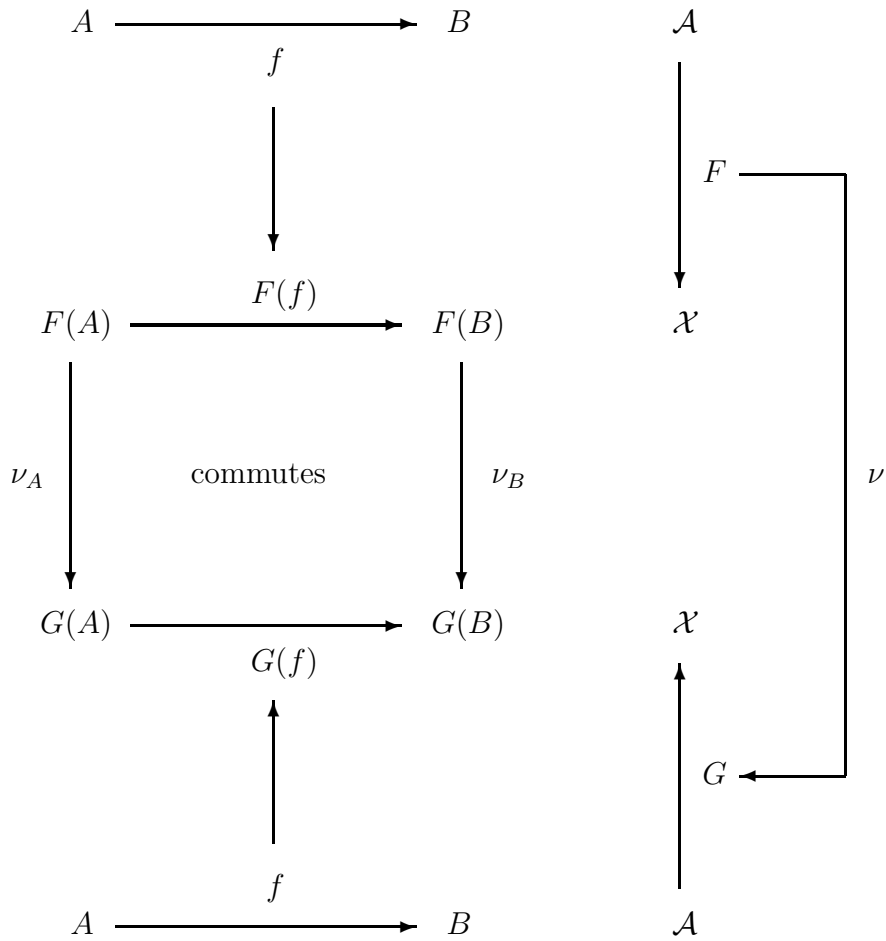
\begin{math}
\setlength{\unitlength}{0.2cm}
\thicklines
\begin{picture}(90,62)

\put(10,60){$A$}
\put(13,60.5){\vector(1,0){20}}
\put(35,60){$B$}
\put(50,60){${\cal A}$}
\put(23,57.5){$f$}
\put(8,40){$F(A)$}
\put(14,40.5){\vector(1,0){18}}
\put(22,42){$F(f)$}
\put(34,40){$F(B)$}
\put(50,40){${\cal X}$}
\put(23.5,55){\vector(0,-1){9.5}}
\put(50.5,58){\vector(0,-1){15}}
\put(52,50){$F$}

\put(10,38){\vector(0,-1){15}}
\put(6,30){$\nu_A$}
\put(36,38){\vector(0,-1){15}}
\put(38,30){$\nu_B$}
\put(18,30){$\mbox{commutes}$}

\put(10,0){$A$}
\put(13,0.5){\vector(1,0){20}}
\put(35,0){$B$}
\put(50,0){${\cal A}$}
\put(23,2.5){$f$}
\put(8,20){$G(A)$}
\put(14,20.5){\vector(1,0){18}}
\put(22,17.5){$G(f)$}
\put(34,20){$G(B)$}
\put(50,20){${\cal X}$}
\put(23.5,6){\vector(0,1){9.5}}
\put(50.5,3){\vector(0,1){15}}
\put(52,10){$G$}

\put(54.5,50.5){\line(1,0){7}}
\put(61.5,50.55){\line(0,-1){40}}
\put(61.5,10.5){\vector(-1,0){7}}
\put(63,30){$\nu$}

\end{picture}
\end{math} \\

We recall that $\nu_A$ and $\nu_B$ above are morphisms in the category ${\cal X}$. \\

A natural transformation $~\nu~$ between categories ${\cal A}$ and ${\cal X}$ is called a {\it
natural isomorphism}, if and only if for every object $A$ in category ${\cal A}$, the
morphism \\

\begin{math}
\setlength{\unitlength}{0.2cm}
\thicklines
\begin{picture}(90,22)

\put(8,20){$F(A)$}
\put(9.6,18){\vector(0,-1){15}}
\put(6,10.5){$\nu_A$}
\put(8,0){$G(A)$}

\end{picture}
\end{math} \\

is an isomorphism in category ${\cal X}$, that is, there exists a morphism $G(A)
~\stackrel{\mu_A}\longrightarrow~ F(A)$ in category ${\cal X}$, such that \\

$~~~~~~~~~~ \mu_A \circ \nu_A ~=~ id_{F(A)},~~~\nu_A \circ \mu_A ~=~ id_{G(A)} $ \\ \\

\end{document}